\documentclass[11pt]{article}

\makeatletter
\@addtoreset{equation}{section}
\makeatother

\newtheorem{theorem}{Theorem}
\newtheorem{algorithm}{Algorithm}
\usepackage{epsfig}
\usepackage{amssymb}
\usepackage{amsmath}
\usepackage{graphicx}

\def\tfrac#1#2{{{\lower.6ex
\hbox{$\scriptstyle#1$}}\over
{\raise.7ex
\hbox{$\scriptstyle#2$}}}}

\def\bigO{{\cal O}}

\def\protectbold#1{\protect{\boldmath{$#1$}}}

\def\Frac#1#2{\frac{\displaystyle{#1}}{\displaystyle{#2}}}

\def\calL{{\cal L}}
\def\bigO{{\cal O}}

\def\Ai{{\rm Ai}}
\def\Bi{{\rm Bi}}
\def\Gi{{\rm Gi}}
\def\Hi{{\rm Hi}}

\def\bigO{{\cal O}}

\begin{document}

\title{Recent software developments for special functions in the 
Santander-Amsterdam project}

\author{Amparo Gil\\
Departamento de Matem\'atica Aplicada y CC. de la Computaci\'on.\\
ETSI Caminos. Universidad de Cantabria. 39005-Santander, Spain.\\
\and
Javier Segura\\
        Departamento de Matem\'aticas, Estad\'{\i}stica y 
        Computaci\'on,\\
        Univ. de Cantabria, 39005 Santander, Spain.\\   
\and
Nico M. Temme\\
IAA, 1391 VD 18, Abcoude, The Netherlands\footnote{
Former address: CWI,
Science Park 123, 1098 XG Amsterdam, The Netherlands}.\\ \\
{ \small
  e-mail: {\tt
    amparo.gil@unican.es,
    javier.segura@unican.es, 
    Nico.Temme@cwi.nl}}
    }

\maketitle
\begin{abstract}
We give an overview of published algorithms by our group and of current activities and future plans. 
In particular, we give details on methods for computing special functions and discuss in detail two 
current lines of research. Firstly, we describe the recent developments for the
computation of central and non-central $\chi$-square cumulative distributions (also called Marcum $Q-$functions), and we present a new quadrature method
for computing them. Secondly, we describe the 
fourth-order methods for computing zeros of special functions recently developed, and we 
provide an explicit example for the computation of complex zeros of Bessel functions. We end with an overview of 
published software by our group for computing special functions. 
\end{abstract}

\vskip 0.8cm \noindent
{\small
Keywords \& Phrases:
special functions; incomplete gamma functions, Marcum's $Q-$function; zeros of special functions; algorithms; numerical 
software.
}

\section{Introduction}\label{sec:intro}
Our project on software developments for special functions started in 1997.
Earlier we published a number of papers in this area, and the plan was to combine all our expertise in order to produce quality software based on a selection of the many methods available for special functions, and by justifying these methods in particular cases with the help of  elements of numerical and mathematical analysis, such as recurrence relations, numerical quadrature and asymptotic expansions.

We discuss current activities in two different lines of research. Firstly, we
discuss the computation of $\chi$-square cumulative distributions and, in particular, we
present a new quadrature method for computing the non-central distribution, also called Marcum's $Q-$function. Secondly, 
we describe in a unified way two recent methods for computing zeros of special functions (for real and complex zeros), we
give an explicit example of computation 
for complex zeros of Bessel functions and discuss plans for software implementations. We end with an overview of 
our published software for computing
special functions. 

First, we  briefly outline the basic methods and principles we consider in the construction of special function
software.

\section{Numerical methods and basic principles}

Our book \cite{Gil:2007:NSF} ``Numerical methods for special functions"
 appeared in 2007 and describes four basic methods that we have used in writing software. These methods are:
convergent and asymptotic series, Chebyshev expansions, linear recurrence relations and quadrature methods. The
book also describes numerical methods for computing continued fractions, methods for computing zeros of special
functions and the computation with uniform asymptotic expansions, among other topics. Usually, several of these methods
have to be combined for computing a special function.

\subsection{Our principles of designing algorithms}\label{sec:princ}
 Our approach in making software can be described by the following principles.  
 
 \begin{enumerate}
 \item
 The main objective is to develop Fortran 90 codes which produce reliable double precision values. 
 We use Maple or Mathematica for obtaining coefficients in expansions and for verifying algorithms, but not usually
 in the final product. The exceptions are some codes for computing zeros of special functions.
	\item
	A given special function is usually a special case of a more general function. 
Our approach is bottom-up and when a simple but important function is demanding software, we prefer to start with the ``simple"
case. For example: Airy functions are special cases of the more general Bessel functions, but we have written codes for the Airy functions themselves.  
       \item We accept that it is necessary to combine several methods in order to compute a function accurately and efficiently
for a wide range of its variables.
	\item We accept that a theoretical error analysis is usually impossible for functions with several real or complex 
variables. We accept more empirical approaches.
	\item The accuracy analysis is usually done by using functional relations, such as Wronskian relations or by comparing
with an alternative method of computation. 
	\item The selection of methods in different parameter domains is based on speed and accuracy, where the latter may 
prevail. For large real or complex parameters  scaling of the result is useful to avoid underflow or overflow in our finite arithmetic environment.
 \end{enumerate}

\section{Incomplete gamma functions, Marcum \protectbold{Q-}function}\label{sec:igm}

As before commented, present interest is focused on distribution functions, in particular on the incomplete gamma
functions (see \cite{Gil:2012:IGR}) and generalizations. Incomplete gamma functions are the 
central $\chi$-square cumulative distributions and are defined by
\begin{equation}\label{eq:int01}
\gamma(a,x)=\int_0^x t^{a-1} e^{-t}\,dt, \quad
\Gamma(a,x)=\int_x^{\infty} t^{a-1} e^{-t}\,dt.
\end{equation}
We concentrate on the ratios
\begin{equation}\label{eq:int02}
P(a,x)=\frac{1}{\Gamma(a)}\gamma(a,x), \quad
Q(a,x)=\frac{1}{\Gamma(a)}\Gamma(a,x),
\end{equation}
where we assume that $a$ and $x$ are positive. We can use the well-known series expansions, asymptotic expansions, recurrence relations, and a method based on the uniform asymptotic representation of these ratios in terms of the complementary error function. Apart from the last method, we mainly use the methods considered in \cite{Gautschi:1979:IGF,Gautschi:1979:CPI}, although Gautschi used a different set of functions and also negative values of $a$. For applications in mathematical statistics and probability theory the ratios in \eqref{eq:int02} are more relevant.
Furthermore, \cite{Gil:2012:IGR} in we described algorithms for inverting the incomplete gamma ratios, 
and the algorithm improves the one given in \cite{Didonato:1986:CIG,Didonato:1987:AFS}.

The results of the algorithms for the incomplete gamma ratios are used in our current project on the the generalized Marcum $Q-$function, which is defined in \eqref{eq:defQmu} and \eqref{eq:defPmu}. 
A paper with full details of our numerical computations has been submitted \cite{Gil:2013:CMQ}.
The relation with the incomplete gamma functions becomes clear when expanding the Bessel function in its power series, which gives
\begin{equation}\label{eq:int04}
Q_\mu(x,y)=e^{-x}\sum_{n=0}^\infty \frac{x^n}{n!} Q(\mu+n,y),
\end{equation}
and
\begin{equation}\label{eq:int05}
P_\mu(x,y)=e^{-x}\sum_{n=0}^\infty \frac{x^n}{n!} P(\mu+n,y).
\end{equation}

\subsection{A new quadrature method for computing the Marcum $Q$-function}\label{sec:quad}
As we have explained in \cite[Chapter~5]{Gil:2007:NSF}, numerical quadrature can be an important tool for evaluating special functions. In particular, when selecting suitable  integral representations for these functions. In Chapter 12 of our book we have given several examples, see also \S\S\ref{sec:kinu}, \ref{sec:pcf}, \ref{sec:leg}. The quadrature methods are usually based on the simple trapezoidal rule, that is very efficient and accurate for certain integral representations that follow from contour integrals in the complex plane, for example through the saddle point of the integrand. We have shown for many cases that the trapezoidal rule is an effective tool for a large domain of the parameters. Usually we use the method to overlap  power series methods near the origin and asymptotic methods for large parameters.

In this section we describe the trapezoidal method for a special case, the Marcum $Q-$function, and we give all analytical details that are needed for including the method in a computer program. In fact, this is generally the way of working: we need quite some analytical preparations before we apply the simple trapezoidal rule, but the extra work is usually very rewarding.

We define the generalized Marcum $Q-$function by using the integral representation
\begin{equation}
\label{eq:defQmu}
Q_{\mu} (x,y)=\displaystyle x^{\frac12 (1-\mu)} \int_y^{+\infty} t^{\frac12 (\mu -1)} e^{-t-x} I_{\mu -1} \left(2\sqrt{xt}\right) \,dt,
\end{equation}
where $x\ge 0$, $y\ge 0$, $\mu>0$ and $I_\mu(z)$ is the modified Bessel function. 
We also use the complementary function
\begin{equation}
\label{eq:defPmu}
P_{\mu} (x,y)=\displaystyle x^{\frac12 (1-\mu)}\int_0^{y} t^{\frac12 (\mu -1)} e^{-t-x} I_{\mu -1} \left(2\sqrt{xt}\right) \,dt,
\end{equation}
and the complementary relation reads
\begin{equation}\label{eq:PQcompl}
P_{\mu}(x,y)+Q_{\mu} (x,y)=1.
\end{equation}

The generalized Marcum $Q-$function, which
for $\mu=1$ reduces to the ordinary Marcum function, is used in problems on radar detection and 
communications; see
\cite{Marcum:1960:AST}.  In this field, $\mu$ is the number of
independent samples of the output of a square-law detector. In
our analysis $\mu$ is not necessarily an integer number. There is an extensive literature regarding this function, also from the areas of statistics and probability theory, where they are called the non central $\chi$-square or the non central gamma distribution.
For papers  concentrating on numerical aspects we refer to 
  \cite{Andras:2011:GMQ,Ashour:1990:OTC,Dyrting:2004:ENC,Helstrom:1992:CGM,Knusel:1996:CNG,Robertson:1976:CON,Ross:1999:ACN,Shnidman:1989:GMF}. As explained earlier in this section, the central gamma distribution is in fact the incomplete gamma function.

The integrals in \eqref{eq:defQmu} and \eqref{eq:defPmu} give stable integral representations, but we prefer a representation in terms of elementary functions. Also, one important point for applying the trapezoidal rule  efficiently is the vanishing of the integrand and many (or all) of its derivatives at the finite endpoint of integration.

In the present case we start with the complex integral representation (see \cite[Eq.~(2.3)]{Temme:1993:ANA})
\begin{equation}\label{eq:quad01}
Q_\mu(x,y)=\frac{e^{-x-y}}{2\pi i}\int_{\calL_Q}  \frac{e^{x/s+ys}}{1-s}\,\frac{ds}{s^\mu},
\end{equation}
where ${\calL_Q} $ is a vertical line that cuts the real axis in a point $s_0$, with $0<s_0<1$. 
For the complementary function we have
\begin{equation}\label{eq:quad02}
P_\mu(x,y)=\frac{e^{-x-y}}{2\pi i}\int_{\calL_P} \frac{e^{x/s+ys}}{s-1}\,\frac{ds}{s^\mu},
\end{equation}
now with  a vertical line ${\calL_P} $ that cuts the real axis at  a point $s_0$ with $s_0>1$. Both representations follow from each other by shifting the contour across the pole at $s=1$ and picking up the residue.

We use scaled variables $x, y$ and write the representation in \eqref{eq:quad01} in the form
\begin{equation}\label{eq:quad03}
Q_\mu(\mu x,\mu y)=\frac{e^{-\mu(x+y)}}{2\pi i}\int_{\calL_Q}  \frac{e^{\mu\phi(s)}}{1-s}\,ds,
\end{equation}
where ${\calL_Q} $ is as in \eqref{eq:quad01}, and
\begin{equation}\label{eq:quad04}
\phi(s)=\frac{x}{s}+ys-\ln s.
\end{equation}

The first step in choosing a suitable contour for numerical quadrature is determining the saddle point of $\phi(s)$ and shifting the contour through this saddle point. We have
\begin{equation}\label{eq:quad05}
\phi^\prime(s)=-\frac{x}{s^2}+y-\frac{1}{s},
\end{equation}
which vanishes at the point  (the negative zero is not relevant here)
\begin{equation}\label{eq:quad06}
s_0=\frac{1+\sqrt{1+\xi^2}}{2y},\quad \xi=2\sqrt{xy}.
\end{equation}
We want to evaluate the $Q-$function when $y > x+1$ (in the scaled variables). When $y<x+1$ we would compute the
$P-$function instead and use (\ref{eq:PQcompl}).

We continue with $y>x+1$ with ${\calL_Q}$ through $s_0$ of \eqref{eq:quad06}.  It is easy to verify that in that case the saddle point $s_0$ 
satisfies $0<s_0<1$. 
The next step is to deform this line into a contour (still through $s_0$) on which $\Im\phi(s)=\Im\phi(s_0)$. This gives the equation
\begin{equation}\label{eq:quad07}
\Im\left(\frac{x}{r}e^{-i\theta}+yre^{i\theta}-\ln r-i\theta\right)=0,
\end{equation}
because $\Im\phi(s_0)=0$. Here we have used polar coordinates
\begin{equation}\label{eq:quad08}
s=r e^{i\theta},\quad r>0,\quad -\pi<\theta<\pi.
\end{equation}
Solving for $r$ we obtain
\begin{equation}\label{eq:quad09}
r(\theta)=\frac{1}{2y}\left(\frac{\theta}{\sin\theta}+\rho(\theta,\xi)\right), \quad \rho(\theta,\xi)=\sqrt{\left(\frac{\theta}{\sin\theta}\right)^2+\xi^2}.
\end{equation}
This defines a parabola-shaped contour through $s_0$ given in \eqref{eq:quad06} and extending to $\Re s\to -\infty$ (when $\theta\to\pm\pi$).

On this contour we have
\begin{equation}\label{eq:quad10}
\phi(s)=\Re \phi(s)=\cos\theta \rho(\theta,\xi)-\ln r(\theta).
\end{equation}
Integrating on the contour with respect to $\theta$ gives
\begin{equation}\label{eq:quad11}
Q_\mu(\mu x,\mu y)=\frac{e^{-\mu(x+y)}}{2\pi i}\int_{-\pi}^\pi  e^{\mu\Re \phi(s)}\frac{r^\prime(\theta)+ir}{e^{-i\theta}-r}\,d\theta,
\end{equation}
and by evaluating the fraction in the integrand we obtain the real representation
\begin{equation}\label{eq:quad12}
Q_\mu(\mu x,\mu y)=\frac{e^{-\mu(x+y)}}{2\pi}\int_{-\pi}^\pi  e^{\mu\Re \phi(s)} f(\theta)\,d\theta, 
\end{equation}
where $f(\theta)$ is an even function (the odd part does not contribute), and is defined by
\begin{equation}\label{eq:quad13}
f(\theta)=\frac{\sin\theta\, r^\prime(\theta)+\left(\cos\theta-r(\theta)\right)r(\theta)}{r^2(\theta)-2r(\theta)\cos\theta +1}.
\end{equation}

This representation has the required form: the integrand vanishes with all its derivatives at the endpoints of the interval. But this is not the final representation. The point is, we like to take out the value $\exp(\mu\Re\phi(s_0))$ at $\theta=0$ in order to avoid numerical instabilities for small values of $\theta$. We write
\begin{equation}\label{eq:quad14}
-x-y+\phi(s_0)=-\tfrac12\zeta^2
\end{equation}
and we obtain our final representation 
\begin{equation}\label{eq:quad15}
Q_\mu(\mu x,\mu y)=\frac{e^{-\frac12\mu\zeta^2}}{2\pi }\int_{-\pi}^\pi  e^{\mu  \psi(\theta)} f(\theta)\,d\theta, 
\end{equation}
where
\begin{equation}\label{eq:quad16}
\psi(\theta)=\cos\theta \rho(\theta,\xi)-\sqrt{1+\xi^2}-\ln\frac{\Frac{\theta}{\sin\theta}+\rho(\theta,\xi)}{1+\sqrt{1+\xi^2}},
\end{equation}
with $ \rho(\theta,\xi)$ defined in \eqref{eq:quad09}.

\subsubsection{Stable computations}\label{sec:quadsub01}
It will be clear that for numerical evaluations we need some preparations for small values of $\theta$. First we observe that
\begin{equation}\label{eq:quad17}
\psi(\theta)=-\tfrac12\sqrt{1+\xi^2}\,\theta^2\left(1+\frac{2-3\xi^2}{36(1+\xi^2)}\theta^2+\bigO\left(\theta^4\right)\right),\quad \theta\to 0,
\end{equation}
and we can compute more terms in this expansion. However, these terms depend on $\xi$, and it is better to make first a few analytical steps. First we write
\begin{equation}\label{eq:quad18}
\cos\theta \rho(\theta,\xi)-\sqrt{1+\xi^2}=\frac{\left(\Frac{\theta}{\sin\theta}\right)^2-1-\theta^2-\xi^2\sin^2\theta}{\cos\theta \rho(\theta,\xi)+\sqrt{1+\xi^2}},
\end{equation}
and for 
\begin{equation}\label{eq:quad19}
\left(\frac{\theta}{\sin\theta}\right)^2-1=\frac{(\theta-\sin\theta)(\theta+\sin\theta)}{\sin^2\theta}
\end{equation}
we may use the simple expansion 
\begin{equation}\label{eq:quad20}
\theta-\sin\theta=\tfrac16\theta^3\left(1+6\sum_{k=1}^\infty (-1)^k\frac{\theta^{2k}}{(2k+3)!}\right).
\end{equation}
Although the convergence is very fast, in the algorithm we use a Chebyshev expansion for this quantity.
 
The logarithmic term in \eqref{eq:quad16} needs also some care. We avoid the straightforward evaluation because when $\theta$ is small the argument tends to unity. In fact, we use a Chebyshev expansion of the function $\ln(1+z)$ that is accurate for small values of $\vert z\vert$. For this we write
\begin{equation}\label{eq:quad21}
z=\frac{\Frac{\theta}{\sin\theta}+\rho(\theta,\xi)}{1+\sqrt{1+\xi^2}}-1=\frac{\Frac{\theta-\sin\theta}{\sin\theta}+
\rho(\theta,\xi)-\sqrt{1+\xi^2}}{1+\sqrt{1+\xi^2}},
\end{equation}
which can be written as
\begin{equation}\label{eq:quad22}
z=\frac{\Frac{\theta-\sin\theta}{\sin\theta}}{1+\sqrt{1+\xi^2}}\left(1+\frac{\Frac{\theta}{\sin\theta}+1}
{\rho(\theta,\xi)+\sqrt{1+\xi^2}}\right).
\end{equation}

Finally, for the function $f(\theta)$ defined in \eqref{eq:quad13} we need a stable representation of $r^\prime(\theta)$. We have from \eqref{eq:quad09}
\begin{equation}\label{eq:quad23}
r^\prime(\theta)=\frac{\sin\theta-\theta\cos\theta}{2y\sin^2\theta}\left(
1+\frac{\theta}{\rho(\theta,\xi)\sin\theta}\right), 
\end{equation}
and we write
$\sin\theta-\theta\cos\theta= \sin\theta-\theta+2\theta\sin^2\left(\frac12\theta\right)$.

\subsubsection{Where to use this quadrature method}\label{sec:wherequad}

The function $f(\theta)$ defined in \eqref{eq:quad13} becomes singular when $s=r(\theta)e^{i\theta}=1$. This happens if $\theta=0$, and otherwise for imaginary values of $\theta$. If $r(0)=1$, then $y=x+1$, and when $\vert y-x-1\vert$ is small we need many function evaluations in the quadrature rule, which we like to avoid. When  $\vert y-x-1\vert$ is small we use asymptotic expansions. A safe domain follows from the asymptotic analysis and is given by
\begin{equation}\label{eq:where01}
x+1-b\sqrt{2/\mu}\,\sqrt{2x+1}<y<x+1+b\sqrt{2/\mu}\,\sqrt{2x+1},
\end{equation}
for some positive number $b$, say $b=1$. Outside this parabolic shaped domain around the line $y=x+1$ we can use for the representation in \eqref{eq:quad15} the trapezoidal rule with great efficiency. 
It should be observed that, although the selection of the contour is based on asymptotic analysis,  the quadrature method does not need necessarily large parameters, although it performs also quite well in that case.

The domain in \eqref{eq:where01} is described in terms of  the scaled variables for $Q_\mu(\mu x,\mu y)$. For the unscaled variables in $Q_\mu(x,y)$ we use
\begin{equation}\label{eq:where02}
x+\mu-b\sqrt{4x+2\mu}<y<x+\mu+b\sqrt{4x+2\mu}.
\end{equation}

In \cite{Helstrom:1992:CGM} the trapezoidal rule is used by including the pole at $s=1$ in \eqref{eq:quad03} in the function $\phi(s)$. That is, by writing (compare \eqref{eq:quad04})
\begin{equation}\label{eq:where03}
\widetilde{\phi}(s)=\frac{x}{s}+ys-\ln s- \frac{1}{\mu}\ln(1-s).
\end{equation}
In that case the saddle point has to be calculated from a cubic polynomial and the  contour follows from  $\Im\widetilde{\phi}(s)=0$ through that saddle point. Helstrom used an approximation of this contour by taking a parabola centered at the saddle point. The  tabled results show correct values at the the critical values $y=x+1$. 

It is of interest to see how the method described in \cite{Chiarella:1968:OTE}, in which the trapezoidal rule is used when a pole is close to the saddle point, can be applied to the integral in \eqref{eq:quad01}.

\section{Zeros of special functions}
\label{sec:zer}

There is ongoing work on software developments for the computation of zeros of special functions. The zeros
of special functions are important quantities appearing in numerous scientific applications. 

The problem we are addressing is, given a (special) function $f(x)$, develop a program for computing with 
certainty the solutions of $f(x)=0$. We consider real zeros, but we also discuss the extension to complex zeros.

We restrict our attention to
functions which are solutions of second-order differential equations 
\begin{equation}
y''(x)+B(x)y'(x)+A(x)y(x)=0.
\end{equation}
Special
cases are, for instance, the classical orthogonal polynomials (for which the zeros are the nodes of Gaussian quadrature formulas), and Bessel functions. 

Except for some particular cases (like, for instance, the Golub-Welsch algorithm for orthogonal
polynomials, see for instance \cite[\S5.3.2]{Gil:2007:NSF}) the development of software for computing 
zeros of a function requires that software for computing the function is available; this is also the case of our algorithms. 
In this sense, the 
computation of zeros of special functions is secondary to the evaluation of the functions themselves.
Most of the
software for the evaluation of special functions is written in Fortran and we consider this programming language
when special function software is available.

Another possibility is using commercial software with built-in special function commands.
Mathematica and Maple are interesting possibilities with a large span of mathematical functions available
We use Maple for some of our packages; current activities include the development of software packages for computing
real or complex zeros of special functions using the fourth-order methods of \cite{Segura:2010:RCZ,Segura:2012:CCZ}.

Next, we describe the basic ingredients for the method introduced in \cite{Segura:2010:RCZ} 
for real zeros and the status of the software implementations, both
in Fortran and in Maple. Here we motivate the method as a direct consequence of
Sturm's comparison theorem \cite[page 19]{Szego:1975:OP}. In the second place, we briefly discuss the extension of the
method for the computation of complex zeros presented in \cite{Segura:2012:CCZ}, we give some specific examples of 
application in Maple and discuss details of the implementation.

\subsection{Computation of real zeros}

The method introduced in \cite{Segura:2010:RCZ} is a fast fourth-order method which is able to compute all the
zeros of any solution of a second-order ODE $y''(x)+A(x)y(x)=0$ \footnote{Any differential equation
$y''(x)+B(x)y'(x)+A(x)y(x)=0$ with $B(x)$ differentiable can be transformed to normal form (no first derivative term) with
the change of function $w(x)=\exp(\frac12 \int^x B(\zeta)d\zeta)y(x)$} in any real interval where $A(x)$
is continuous, provided the monotonicity properties of $A(x)$ in this interval are known in advance. This method
can be motivated as a direct consequence of the following Sturm theorem:

\begin{theorem}[Sturm comparison]
Let $y(x)$ and $w(x)$ be solutions of $y''(x)+A_{y}(x)y(x)=0$ and $w''(x)+A_{w}(x)w(x)=0$ respectively, with
 $A_w (x)>A_{y}(x)>0$. If $y(x_0)w'(x_0)-y'(x_0)w(x_0)=0$ and $x_y$ and $x_w$ are the zeros of $y(x)$ and $w(x)$
closest to $x_0$ and larger (or smaller) than $x_0$, then $x_w<x_y$ (or $x_w>x_y$).
\end{theorem}

\begin{figure}[tb]
\vspace*{1cm}
\begin{center}
\centerline{\protect\hbox{\psfig{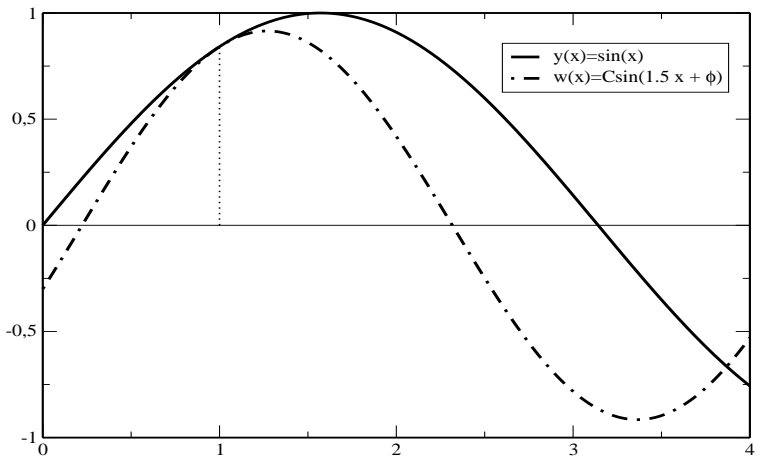}}}
\end{center}
\caption{The function $y(x)=\sin (x)$ is a solution of $y''(x)+A_y(x) y(x)=0$, $A_y(x)=1$, and 
$w(x)=C\sin(1.5\,x+\phi)$, $C=0.9153\ldots$, $\phi=-0.3336\ldots$
is a solution of  $w''(x)+A_w(x) w(x)=0$, $A_w(x)=2.25$. The curves are tangent at $x_0=1$, where
$y(x_0)w'(x_0)-y'(x_0)w(x_0)=0$. Because $A_w (x)>A_{y}(x)>0$, the zeros of $w(x)$ are closer to
$x_0=1$ than the zeros of $y(x)$.
}
\label{fig1}
\end{figure}

The proof is straightforward. An intuitive explanation is as follows. 
The differential equations of the form $y''(x)+A(x)y(x)=0$ have
solutions which may oscillate if $A(x)>0$, and the oscillations are more rapid as $A(x)$ is larger. In the theorem,
because $A_w(x)>A_{y}(x)>0$, the solutions of the second equation oscillate more rapidly and their zeros tend to be
closer together.  Now, because we have the hypothesis 
$y(x_0)w'(x_0)-y'(x_0)w(x_0 )=0$ we can consider that $y(x_0 )=w(x_0)$ and $y'(x_0 )=w'(x_0)$, and
there is no loss of generality
because a solution of a linear homogeneous ODE can be multiplied by a constant and it remains a solution of the same ODE
with the same zeros.
Therefore, the solutions $y(x)$ and $w(x)$ have the same initial conditions at $x_0$; but because $w(x)$ oscillates more rapidly 
than $y(x)$, then then there is necessarily a zero of $w(x)$ between $x_0$ and any zero of $y(x)$.
 This is illustrated graphically in Figure \ref{fig1}, where we consider the simple case of equations with
constant coefficients.

From Sturm's theorem we can construct a method for the computation of the zeros
of solutions of $y''(x)+A(x)y(x)=0$ when $A(x)$ is monotonic. If $A(x)$ is a decreasing (increasing)
function and $A(x)>0$ we compute the zeros with an increasing (decreasing) sequence; if $A(x)<0$ the solutions have one zero
 at most. Given a value $x_0$, the zero of $y(x)$ closest to $x_0$ and larger (smaller) than $x_0$ can be computed with 
certainty using the following scheme.

\begin{algorithm}[Zeros of $y''(x)+A(x)y(x)=0$, $A(x)>0$ monotonic].
\label{Algo1}

Let $x_0<\alpha$ with $y(\alpha )=0$ and such that there is no zero of $y(x)$ 
between $x_0$ and $\alpha$.
  
Assume that $A(x)$ is decreasing (increasing). Starting from $x_0$, compute $x_{n+1}$ from $x_n$ as follows: 
find a non-trivial solution of the equation $w''(x)+A(x_n)w(x)=0$ such that
$y(x_n)w'(x_n)-y'(x_n)w(x_{n})=0$. Take as $x_{n+1}$ 
the zero of $w(x)$ closest to $x_n$ and larger (smaller) than $x_n$. Then, the sequence $\{x_n\}$ converges to $\alpha$.

\end{algorithm}

This is a direct consequence of Sturm comparison. Consider, for instance, the case of $A(x)$ increasing;
$A(x_n)>A(x)$ for $x>x_n$ and the solutions of $w''(x)+A(x_n)w(x)=0$ 
oscillate more rapidly than those of $y''(x)+A(x)y(x)=0$, and because $y(x_n)w'(x_n)-y'(x_n)w(x_{n})=0$ we
have $x_n<x_{n+1}<\alpha$. For increasing $A(x)$ the iterations with 
the algorithm produce an increasing sequence with upper bound the zero that is computed;
therefore it converges, and necessarily to the zero because, as we will see,
the only fixed point of the resulting iteration function (Eq. (\ref{def1})) in $(x_n,\alpha]$ is precisely $\alpha$. 
The situation is similar to Figure \ref{fig1} ($x>1$), where $y(x)$ would be the special
function and $w(x)$ would provide the next iteration.

 Of course, the algorithm can be applied successive times to generate a sequence of zeros.
\begin{algorithm}[Computing a sequence of zeros, $A(x)$ monotonic].
\label{Algo2}

Let $\alpha_1, \alpha_2$ be consecutive zeros of $y(x)$, with $\alpha_1 <\alpha_2$. 

If $A(x)$ is decreasing and $\alpha_1$ is known, the zero $\alpha_2$ can be computed using
Algorithm \ref{Algo1} with starting value $x_0=\alpha_1$ 
(the first iteration being $x_1=\alpha_1+\pi/\sqrt{A(\alpha_1)}$).

 If $A(x)$ is increasing and $\alpha_2$ is known, the zero $\alpha_1$ can be computed using
Algorithm \ref{Algo1} with starting value $x_0=\alpha_2$ 
(the first iteration being $x_1=\alpha_2-\pi/\sqrt{A(\alpha_2)}$).
\end{algorithm}

As commented before the sequences generated are increasing (decreasing) if $A(x)$ is decreasing (increasing).

The iteration of Algorithm \ref{Algo1} can be explicitly written as follows:

\begin{equation}
\label{def1}
T(x)=x-\Frac{1}{\sqrt{A(x)}}\arctan_{j}(\sqrt{A(x)}h(x)),
\end{equation}
with $h(x)=y(x)/y'(x)$, $j=\mbox{sign}(A'(x))$ and
\begin{equation}
\label{def2}
\arctan_{j}(\zeta)=\left\{
\begin{array}{l}
\arctan(\zeta)
\mbox{ if }j\zeta > 0,\\ 
\arctan(\zeta)+j\pi \mbox{ if }j\zeta\le 0,\\ 
j\pi/2 \mbox{ if } \zeta=\pm \infty.
\end{array}
\right.
\end{equation}
Observe that the only fixed points of $T(x)$ are the zeros of $y(x)$.

The method is of order four which means, roughly speaking, that the number
of correct digits is multiplied by four in each iteration provided the iterations are close enough to the zero. On the
other hand, it has good non-local behavior (see \cite[Definition 4.1]{Segura:2010:RCZ}), 
and once a zero is computed, few iterations are generally
needed to reach a close estimation of the next zero. The fixed point method depends on an arctangent function but 
defined with a different range depending
on whether $A(x)$ is decreasing or increasing. When an iteration close to the zero that is being computed 
is reached, it is better to switch 
to the standard arctangent functions (range $(-\pi /2, \pi /2)$), which gives a fixed point method continuous at the zero
and convergent around the zero

The algorithms need some a priori analysis: the monotonicity properties of the coefficient $A(x)$ must be known in advance, because
the method has to be applied separately in those subintervals where $A(x)$ is monotonic.
This analysis has been completed for
hypergeometric functions \cite{Deano:2004:NIC}. 

  An excellent benchmark for these methods is Maple, which permits constructing algorithms for computing zeros of
a large number of special functions; in \cite{Segura:2010:RCZ}, details on the performance of these methods are given 
for a good number of special functions. It is shown that only with three or four iterations per root, one can attain
one hundred digits accuracy. Following these earlier tests, we are currently developing a Maple package for computing
zeros of special functions \cite{Gil:2012:MPC}. By using the method that we are describing in the next section, this package will later be extended to complex zeros.

  Recently, we developed Fortran programs for computing zeros of Bessel functions ${\cal C}_{\nu}(x)$ 
(solutions of $x^2 y''(x)+xy'(x)+(x^2 -\nu^2)y(x)=0$), of their first derivative
and of the combination  $x{\cal C}_{\nu}+\gamma {\cal C}^{\prime}_{\nu}(x)$. In these methods, we considered
the Fortran 77 codes of Amos \cite{Amos:1986:PPB} for computing Bessel functions. An update to Fortran 90
will be considered, together with new Fortran algorithms for computing zeros of other functions for which 
Fortran programs are available, like parabolic cylinder functions \cite{Gil:2006:CRP,Gil:2011:PCW}, 
Legendre functions \cite{Burry:2010:NF9}, conical functions \cite{Gil:2012:AIA}, and 
modified Bessel functions of imaginary argument~\cite{Gil:2004:AMB}.

\subsection{Complex zeros}

 It is possible to extend the previous fourth-order method to zeros in the complex plane, although it is not easy
to prove its convergence in full generality. However, the WKB approximation (also named Liouville approximation)
\cite[Chap.~6]{Olver:1997:ASF} motivates why the method for complex zeros works \cite{Segura:2012:CCZ}.

Let us start by considering the trivial case of $A(z)$ constant. Then the
general solution of $y''(z)+A(z)y(z)=0$ reads
$$
y(z)=C\sin \left(\sqrt{A(z)}\,(z-\psi)\right),
$$
and the zeros are over the line
$$
z=\psi+e^{-i\frac{\varphi}{2}}\lambda,\ \lambda\in{\mathbb R},\  \varphi=\arg{A(z)}.
$$
In other words, writing $z=u+iv$ we have that the zeros are over an integral line of
\begin{equation}
\label{LGA}
\Frac{dv}{du}=-\tan (\varphi/2 ).
\end{equation}

Of course, in general $A(z)$ will not be a constant.
The method for complex zeros is based on the assumption that the curves where the zeros lie 
are also given by (\ref{LGA}), but with variable $\varphi$. This assumption is equivalent to considering
that the WKB approximation is accurate. The WKB approximation with a zero at $z^{(0)}$ is
$$
y(z)\approx C A(z)^{-1/4}\sin\left(\int_{z^{(0)}}^z A(\zeta)^{1/2}d\zeta \right).
$$
Then, if $z^{(0)}$ is a zero, other zeros lie over the curve such that
\begin{equation}
\label{ASL}
\Im \int_{z^{(0)}}^z A(\zeta)^{1/2}d\zeta=0,
\end{equation}
and those curves are also given by (\ref{LGA}). These are the so-called anti-Stokes lines
(ASLs).

The method for computing complex zeros precisely follows the path of the ASLs
and it is similar to the method for real zeros.
Given $z^{(0)}$ ($y(z^{(0)})=0$) and assuming that $|A(z)|$ decreases
for increasing $\Re z$, we consider the following algorithm to compute
the next zero $z^{(1)}$:

\begin{algorithm}[Basic algorithm for complex zeros; $|A(z)|$ decreasing].

\begin{enumerate}
\label{algostep}
\item{}Take $z_0=H^+ (z^{(0)})=z^{(0)}+\pi/\sqrt{A(z^{(0)})}$.
\item{} Iterate $z_{n+1}=T(z_n)$ until $|z_{n+1}-z_n|<\epsilon$, with
\begin{equation}
\label{tzeta}
T(z)=z-\Frac{1}{\sqrt{A(z)}}\arctan \left(\sqrt{A(z)}\Frac{y(z)}{y'(z)}\right).
\end{equation}
\item{Take} as approximate zero $z^{(1)}=z_{n+1}$. 
\end{enumerate}
\end{algorithm}

The algorithm can be repeated in order to compute subsequent zeros.
If $|A(z)|$ is decreasing, the same algorithm can be considered but with the first line
replaced by $z_0=H^- (z^{(0)})=z^{(0)}-\pi/\sqrt{A(z^{(0)})}$.

Observations:
\begin{enumerate}
\item{}If $A(z)$ has slow variation, the first step could be 
a good approximation to $z^{(1)}$. In addition, the step is tangent to the anti-Stokes
line (ASL) at $z^{(0)}$ (that is: the straight line joining $z^{(0)}$ and $z_0$ is tangent at $z^{(0)}$ 
to the ASL passing through this point).

\item{}$|z_0-z^{(0)}|<L$ where $L$ is the length of the anti-Stokes arc
between $z^{(0)}$ and the next zero. It is a step by defect, and in the correct direction. 

\item{}$T(z)$ is a fixed point iteration with order of convergence $4$. This fact does not depend on
the validity of the WKB approximation.
\end{enumerate}

\begin{figure}[tb]
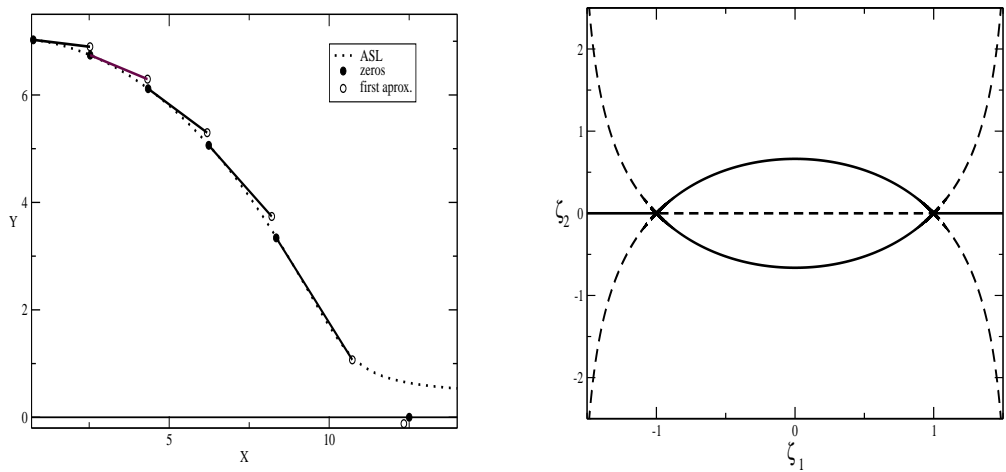

\vspace*{0.8cm}
\begin{center}
\begin{minipage}{4cm}
\centerline{\protect\hbox{\psfig{file=pasos.eps,angle=0,height=6cm,width=6cm}}}
\end{minipage}
\hspace*{3cm}
\begin{minipage}{4cm}
\centerline{\protect\hbox{\psfig{file=stokesbes.eps,angle=0,height=6.2cm,width=6cm}}}
\end{minipage}
\end{center}
\caption{{\bf Left:} Zeros of the Bessel function $Y_{10.35}(z)$ in the quadrant $\Re z>0$, $\Im z>0$ (black circles),
and first estimations to the zeros (white circles); the dotted line is the anti-Stokes line passing through the
zero with larger imaginary part. {\bf Right:} Principal Stokes (dashed lines) and anti-Stokes lines (solid lines) for the 
scaled Bessel equation $d^2 y/d\zeta^2 +\lambda^2 (1-\zeta^{-2})y=0$, $\lambda=\sqrt{\nu^2-1/2}$, $|\nu|>1/2$.
The lines for Bessel functions of order $|\nu|>1/2$ 
correspond to the variable $z=\lambda\zeta$, $\lambda=\sqrt{\nu^2-1/2}$.
}
\label{Fig2}
\end{figure}

Figure \ref{Fig2}, left, shows the complex zeros of the Bessel function
$Y_{10.35}(z)$ in the first quadrant, the first estimations provided 
by the method together with the ASL passing through the zero
with largest imaginary part. The algorithm starts with this zero and 
computes the following zeros (with successively smaller imaginary parts). 
The zeros are very close to
the ASL and the first estimations are very reasonable
with one exception: after computing the last zero with positive imaginary
part, the estimation for the next zero (which is on the real line) 
is not accurate and, furthermore, this zero appears well separated from the ASL. 
We conclude that the WKB approximation
works initially well, but that it is not accurate for computing the last zero 
(although the iteration finally converges to this zero). The problem with this last 
zero is that WKB fails
as a principal Stokes line is crossed. We next explain the notion of principal lines.

A Stokes line (SL) passing through $z_0$ for a differential equation $y''(z)+A(z)y(z)=0$
is given by 
\begin{equation}
\label{SL}
\Re \int_{z_0}^z A(\zeta)^{1/2}d\zeta=0 .
\end{equation}
Compare this with the definition of ASLs (\ref{ASL}). 

Given a point $z_0$ in the complex plane such that $A(z_0)\neq 0,\infty$, then
 there is one and only one ASL (or SL) passing through $z_0$. The situation is 
different if $z_0$ is a zero or a singularity; in particular, if $z_0$ is a zero
of $A(z)$ with multiplicity $m$, there exist $m+2$ ASLs (and SLs) emerging from
$z_0$. We call the ASLs (or SLs) emerging
from the zeros of $A(z)$ principal ASLs (SLs).

For the case of real zeros described in the previous section, the real interval
where the zeros are sought has to be divided in different subintervals where
$A(x)$ is monotonic and the direction of computation is chosen in such a way that
this coefficient decreases. For complex zeros, a similar procedure has to be implemented,
where the position with respect to the principal lines has to be analyzed. Summarizing,
the following strategy proves to be reliable:

The strategy combines the use of 
$H^{(\pm)}=z\pm \Frac{\pi}{\sqrt{A(z)}}$ and $T(z)$ (Eq.~(\ref{tzeta}))
following these rules:

\begin{enumerate}
\item{Divide} the complex plane into disjoint domains separated by the principal ASLs and SLs and compute separately in each domain. 
\item{In each domain}, start away from the principal SLs, close to a principal ASL and/or 
singularity (if any). Iterate $T(z)$ until a first zero is found. If a value outside the domain is reached, stop the search in that domain. 
\item{}Proceed with the basic algorithm, choosing the displacements $H^{(\pm)}(z)$ in the direction of approach
to the principal SLs and/or singularity. 
\item{Stop} when a value outside the domain is reached.
\end{enumerate}

No exception has been found (so far tested for Bessel functions, parabolic cylinder functions and Bessel polynomials).

As an illustration, let us consider the computation of the complex zeros of Bessel functions, and particularly, of
the  Bessel function of the second kind $Y_{\nu}(z)$. The principal Stokes and ASLs are shown in Fig. \ref{Fig2}, 
right. Because the zeros follow ASLs, we expect that the zeros could lie over the principal ASLs 
(the eye-shaped region and part of the real axis) or inside each of the domains separated by these lines. 
The eye-shaped region cuts the imaginary axis at $\zeta=i c$ where $c$ are the real roots of 
$s-\frac12\log((s+1)/(s-1)=0$, $s=\sqrt{1+c^2}$, that is $c=\pm 0.66274321...$. In the particular case of $Y_{\nu}(z)$,
$|\nu| >1/2$, there are zeros very close to the eye-shaped region, in part of the positive real axis and close to the negative real
axis.

First we discuss the computation of the zeros of $Y_{\nu}(z)$ for real orders $\nu>1/2$; for other Bessel 
functions and real orders
the algorithm will be essentially the same. For computing the zeros in the region satisfying
 $\Im z\ge 0$ and $|\Re  z |<L$, with $L$ a given positive number, 
we can consider the following steps:
\begin{enumerate}
\item{Zeros on the positive real line or above the line}: 
take the starting value $z=L+i$, and compute a first zero by iterating $T(z)$. Then, 
proceed using Algorithm \ref{algostep} with 
$H^{-}$ until a value $z$ with 
$\Re z <\sqrt{\nu^2-1/4}$ or $\Im z <-\epsilon$, with $\epsilon$ a small positive number, 
is reached.\footnote{The initial value $z=L$ also works in this case,
but for functions different from $Y_{\nu}(z)$ it is preferable to initiate with $z=L+id$, $d$ a positive number, in order to compute zeros above the positive
real line.}
\item{Zeros above the negative real axis:} take the starting value  $z=-L+i$ and proceed similarly as for the positive real zeros, but with $H^+$.
\item{Zeros on, above or inside the eye-shaped region:} 
start with the estimation $z=0.663\sqrt{\nu^2-1/4}$, compute a first zero 
and move to the right with $H^+$ until a value $z$ such that 
$\Re z >\sqrt{\nu^2-1/4}$ is reached. Similarly, with the same starting zero, move to the left with $H^-$ until a value until 
$\Re z <-\sqrt{\nu^2-1/4}$ is reached. It is also convenient to have, as for the real zeros, the stopping criterion  $\Im z <-\epsilon$.
\end{enumerate}

For the zeros with $\Im z<0$, the same strategy as for $\Im z>0$ can be used, mutatis mutandis. This strategy, in fact, works for
any Bessel function of real order, and not only for $Y_{\nu}(z)$.

As an explicit example, we give some explicit Maple instructions to compute some of these zeros. Let us recall that the coefficient for Bessel functions is $A(z)=1-(\nu^2-1/4)/z^2$, for $y(z)=z^{1/2}Y_{\nu}(z)$ (or $y(z)=z^{1/2}J_{\nu}(z)$ and any 
linear combination). For computing the zeros 
of $Y_{\nu}(z)$ over the eye-shaped region in the first quadrant, we consider the following set of instructions.
{\small
\begin{verbatim}
 
     %% Definition of function and coefficients
     %% Parameters: NI, number of iterations;
     %% NI:=3 gives more than 20 correct digits 
     %% a, order of the Bessel function
> restart:NI:=3:Digits:=30:f:=sqrt(z)*BesselY(a,z):
> h:=f/diff(f,z):coef:=sqrt(1-(a^2-0.25)/z^2):
     

     %% Definition of the iteration T and the displacements H+ and H-
> T:=z-1/coef*arctan(coef*h):Hplus:=z+Pi/coef:Hminus:=z-Pi/coef:a:=10.35:

     %% Computation of a first real zero close to x=40
> x:=40:
> for i from 1 to NI do;
>     x:=evalf(subs(z=x,T));
> end do:
> k:=0:


     %% Computing the real zeros in decreasing order.
     %% The algorithm stops when abs(Re(x))<ll
> ll:=evalf(sqrt(a^2-1/4)):
> while abs(Re(x))>ll do;
>     xc[k+1]:=x;k:=k+1;x:=evalf(subs(z=xc[k],Hminus));
>     for i from 1 to NI do;
>          x:=evalf(subs(z=x,T));
>     end do;
> end do:

     %% Computing the first zero close to the eye-shaped region
> x:=0.662*sqrt(a^2-0.25)*I:
> for i from 1 to NI do;
>     x:=evalf(subs(z=x,T));
> end do:    

     %% Computing the rest of zeros close to the eye-shaped region
     %% The algorithm stops when abs(Re(x))<ll
> while Re(x)<abs(ll) do;
>     xc[k+1]:=x;k:=k+1;x:=evalf(subs(z=xc[k],Hplus));
>     for i from 1 to NI do;
>          x:=evalf(subs(z=x,T));
>     end do;
> end do:

     %% Values of the first 9 real zeros and zeros on the first quadrant
> for k from 1 to k do;xc[k];end do;
                       40.8614543431052914942278668204
                       37.6046889871138442191648203069
                       34.3245525472278464463325220636
                       31.0125895901705525894688274748
                       27.6553314757947377301083091987
                       24.2295152087835610891741598170
                       20.6897425044009133146054974078
                       16.9265377476005389175490477986
                       12.5006643034017892734448978362
                       
    0.751714164454579882038334168282 + 7.02610400031859555166912384150 I
    2.52972140044497327028369512294 + 6.74083920492317316987229531276 I
    4.34222349512990658212287934610 + 6.11572064771947423804814026786 I
    6.23553025542440261192777888137 + 5.06294427710806101928204374876 I
    8.34414329510400851874707241949 + 3.33927403838163677352156909007 I
\end{verbatim}
}

\section{Published algorithms and related papers}\label{sec:algo}
In this section we give an overview of our published algorithms and we mention  the papers that explain in detail the methods of computation for these codes.

A combination of numerical methods is needed for computing function values. Apart from the use of convergent and asymptotic
series, the most frequent methods are those based on the use of linear recurrence relations and of numerical
quadrature. 

The classical reference for the computation using three-term recurrence relations is
\cite{Gautschi:1967:CAT} and we also describe different computation methods in \cite[Chapter~4]{Gil:2007:NSF}. In the
use of recurrences, it is crucial to study the conditioning of the recursion (depending on the direction) and we have 
studied this problem in a series of papers, both in the case of Gauss hypergeometric functions and in the Kummer case. See
\cite{Deano:2008:IMD,Deano:2010:CPT,Gil:2006:ABC,Gil:2006:NSS,Segura:2008:NSS}.
  		
Quadrature methods have been used for  Airy and Scorer functions, certain Bessel and Legendre functions, and parabolic cylinder functions. The trapezoidal rule is particularly efficient for computing numerically many special function integral representations, particularly those arising from saddle point methods. In  \S\ref{sec:quad} we have given an example for a finite interval, and we have given more details in \cite{Gil:2003:CSF},  \cite[Chapter~5]{Gil:2007:NSF}, and in earlier cited papers. 

\subsection{Airy and related functions}\label{sec:airy}
We started this topic with the related functions, also called inhomogeneous Airy functions or Scorer functions, which are solutions of the differential equations
\begin{equation}\label{eq:airy01}
\frac{d^2}{dz^2}w(z)-z\,w(z)=\pm\frac{1}{\pi}.
\end{equation}
With the $+$ sign the standard solution is denoted by $\Gi(z)$, with $-$ sign by $\Hi(z)$. Standard solutions of the homogeneous equation are denoted by $\Ai(z)$ and $\Bi(z)$. For details on these functions we refer to \cite{Olver:2010:ARF}. 

For the Scorer functions we have described a number of contour integrals in the complex plane \cite{Gil:2001:ONI} and the corresponding algorithms can be found in \cite{Gil:2002:AGH}.
For the Airy functions we have also used quadrature methods, see for the analysis \cite{Gil:2002:ANQ} and for the algorithms \cite{Gil:2002:AIZ}. We have also investigated the zeros of the Scorer functions \cite{Gil:2003:OZS} and we have derived asymptotic expansions of these zeros.

In the same period other publications appeared on using asymptotic expansions of the Airy functions and the differential equation \cite{Fabijonas:2004:CCA} with algorithms in \cite{Fabijonas:2004:ALG}.

\subsection{Modified Bessel functions of imaginary order}\label{sec:kinu}
The modified Bessel functions of the third kind of purely
imaginary order $K_{ia}(x)$ is used in a number of problems from physics, and it is also the kernel of the Kantorovich-Lebedev transform. The function $I_{ia}(x)$ is not real, but the function
\begin{equation}
L_{ia}(x)=\tfrac{1}{2}\left[I_{-ia}(x)+I_{ia}(x)\right],
\label{definition}
\end{equation}
is a real valued numerically satisfactory companion to $K_{ia}(x)$ in
the sense considered
in \cite[pp. 154--155]{Olver:1997:ASF}.
The Wronskian relation for these functions is  ${{\cal W}}\left[K_{ia}(x),L_{ia}(x)\right]=1/x$.

In \cite{Gil:2003:CMB,Gil:2004:CSM} we have described the analytical details of the algorithms, and the codes are given in \cite{Gil:2004:AMB}. We have used power series representations for small values of $x$, asymptotic expansions for large $x$, Airy-type asymptotic expansions for large $x$ near the turning point $x=a$, and numerical quadrature. Several non-oscillating integral representations have been used that can be obtained from contour integrals and by using saddle point methods.

\subsection{Parabolic cylinder functions}\label{sec:pcf}
These are solutions of the differential equations
\begin{equation}
\label{pcf}
\frac{d^2}{dx^2}w(x)-\left(\pm \tfrac14x^2 +a\right)w(x)=0.
\end{equation}
The equation with the plus sign, with $U(a,x)$, $V(a,x)$ as two independent solutions, has for
$a<0$ two turning points at $z=\pm\sqrt{-2a}$; for large negative values of $a$ uniform asymptotic representations in terms of Airy functions are available. Oscillations occur between the turning points. Compare this with the Hermite polynomial case when $a=-n-\frac12$ with $n=0,1,2,\ldots$. We have investigated many stable integral representations of these functions in \cite{Gil:2004:IRC} for using numerical quadrature. For the corresponding algorithms for $U(a,x)$, $V(a,x)$ and  their derivatives for real parameters we refer to 
\cite{Gil:2006:ARP,Gil:2006:CRP}.

In \cite{Gil:2011:FAC} we have considered the functions $W(a,x)$ and $W(a,-x)$, which are two linearly independent real solutions of the differential equation (\ref{pcf}) with the minus sign.
In this case the oscillations occur outside the  interval $[-\sqrt{2a},\sqrt{2a}]$.
In \cite{Gil:2011:PCW} we have given the algorithms for computing the $W-$functions and their derivatives. 

In the algorithms for the parabolic cylinder functions we have used recursion,  quadrature, and series expansions, including Maclaurin, local Taylor, Chebyshev and Airy-type asymptotic expansions. By factoring the dominant exponential factor, scaled functions could be computed to avoid overflow/underflow limitations. In this way rather large parameter domains in the $(a,x)$ plane could be covered.

\subsection{Legendre functions: toroidal and conical}\label{sec:leg}
The toroidal functions  are $P^m_{n-\frac12}(x)$ and $Q^m_{n-\frac12}(x)$ and appear in the 
solution of Dirichlet problems with toroidal symmetry. For these functions we have used recurrence relations. For the backward recursions we have used starting values from uniform asymptotic expansions valid for fixed $n$ and large $m$; the expansion is uniformly valid for large positive $x$. For the codes and related papers we refer to \cite{Gil:1998:CEP,Gil:2000:ETH,Gil:2001:DTO}.
 
The conical function $P^ \mu_{-1/2+i\tau}(x)$ is also an element of the class of associated Legendre functions, and the combination of the parameters makes it real for real values of $x$.
This function is the kernel of the Mehler-Fock transform, which has numerous applications, and this function is also used in quantum physics, in particular describing the amplitude for Yukawa potential 
scattering. We have used recurrence relations 
with respect to $\mu=m$ (integer) and related continued fractions, and we have discussed the use of forward and backward recursions.  For $x\in(-1,1)$ we have used quadrature methods and uniform expansions for large $\mu$ in terms of elementary functions. For $x\ge1$ we have used uniform expansions in terms of modified $K-$Bessel functions with purely imaginary order, in particular for describing the behavior of $P^ \mu_{-1/2+i\tau}(x)$ near the turning point $x=\sqrt{\tau^2+\mu^2}/\tau$. The methods of computation and  the algorithms can be found in
\cite{Gil:2009:CCF,Gil:2012:AIA}.

\section*{Acknowledgements}
The authors thank the anonymous referees for helpful comments.
The authors acknowledge support from  {\emph{Ministerio de Econom\'{\i}a y Competitividad}, Spain}, 
projects MTM2009-11686 and MTM2012-34787. 

\bibliographystyle{plain}
\bibliography{gstscp}
\end{document}